\magnification=1200
\overfullrule=0pt
\centerline {\bf A note on the Neumann problem}\par
\bigskip
\bigskip
\centerline {BIAGIO RICCERI}\par
\bigskip
\bigskip
\bigskip
\bigskip
Let $\Omega\subset {\bf R}^n$ be a bounded domain with a boundary of class
$C^1$,
let $\alpha:\Omega\to {\bf R}$ be a bounded measurable function with $\inf_{\Omega}\alpha>0$ 
and let $1<p\leq n$. On the Sobolev space $W^{1,p}(\Omega)$, we
consider the norm
$$\|u\|=\left ( \int_{\Omega}|\nabla u(x)|^p dx+
\int_{\Omega}\alpha(x)|u(x)|^pdx\right ) ^{1\over p}\ .$$
 We denote by ${\cal A}_p$ the class of all
Carath\'eodory functions $\varphi:\Omega\times {\bf R}\to {\bf R}$ such that
$$\sup_{(x,\xi)\in \Omega\times {\bf R}}{{|\varphi(x,\xi)|}\over
{1+|\xi|^q}}<+\infty\ ,$$
for some $q$ with  $0<q< {{pn-n+p}\over {n-p}}$ if $p<n$ and $0<q<+\infty$ if $p=n$.
Given $\varphi\in {\cal A}_p$, consider the following Neumann problem
$$\cases {-\hbox {\rm div}(|\nabla u|^{p-2}\nabla u)+\alpha(x)|u|^{p-2}u=
\varphi(x,u)
 & in
$\Omega$\cr & \cr {{\partial u}\over {\partial\nu}}=0 & on
$\partial \Omega$\cr}\eqno{(P_{\varphi})} $$
where $\nu$ denotes the outward unit normal to $\partial\Omega$.\par
\smallskip
 Let us recall
that a weak solution
of $(P_{\varphi})$ is any $u\in W^{1,p}$ such that
 $$\int_{\Omega}|\nabla u(x)|^{p-2}\nabla u(x)\nabla v(x)dx
+\int_{\Omega}\alpha(x)|u(x)|^{p-2}u(x)v(x)dx
-\int_{\Omega}\varphi(x,u(x))v(x)dx=0$$
for all $v\in W^{1,p}(\Omega)$.\par
\smallskip
The functionals $\Phi, J_{\varphi}:W^{1,p}(\Omega)\to {\bf R}$ defined by
$$\Phi(u)={{1}\over {p}}\|u\|^p$$
$$J_{\varphi}(u)=
\int_{\Omega}\left ( \int_0^{u(x)}\varphi(x,\xi)d\xi\right ) dx$$
are $C^1$ in $W^{1,p}(\Omega)$, with derivatives given by
$$\Phi'(u)(v)=\int_{\Omega}|\nabla u(x)|^{p-2}\nabla u(x)\nabla v(x)dx$$
$$J'_{\varphi}(u)(v)=\int_{\Omega}\varphi(x,u(x))v(x)dx$$
for all $u,v\in W^{1,p}(\Omega)$. Consequently, the weak solutions of problem $(P_{\varphi})$ are exactly the critical points in
$W^{1,p}(\Omega)$ of the functional $\Phi-J_{\varphi}$. Moreover,
$J'_{\varphi}$ is compact, while $\Phi'$ is a homeomorphism
between $W^{1,p}(\Omega)$ and its dual.\par
\smallskip
In recent years, the general abstract results established in
[30], [31] (see also [35] for a revisited version) and [33] have widely been applied to study existence and multiplicity of solutions for problems like $(P_{\varphi})$. In this connection, we refer, for instance, to [1]-[29], [32], [34], [37].\par
\smallskip
Very recently, in [36], we obtained a new general three critical points theorem. Before recalling it, let us give the following definition.\par
\smallskip
 If $X$ is a real Banach space, we
denote by ${\cal W}_X$ the class of all functionals $\Phi:X\to {\bf R}$
possessing the following property: if $\{x_n\}$ is a sequence
in $X$ converging weakly to $x\in X$ and
$\liminf_{n\to \infty}\Phi(x_n)\leq \Phi(x)$, then $\{x_n\}$ has a subsequence converging strongly to $x$.\par
\smallskip
Also, given an operator $S:X\to X^*$, 
we say that $S$ admits a continuous inverse on
$X^*$ if there exists a continuous operator $T:X^*\to X$ such that
$T(S(x))=x$ for all $x\in X$.\par
\medskip
THEOREM A ([36], Theorem 1). - {\it Let $X$ be a separable and reflexive real Banach space;
 $\Phi:X\to {\bf R}$ a coercive, sequentially
weakly lower semicontinuous $C^1$ functional, belonging to ${\cal W}_X$,
bounded
on each bounded subset of $X$ and whose derivative admits a continuous inverse
on $X^*$; $J:X\to {\bf R}$ a $C^1$ functional with compact derivative.
Assume that $\Phi$ has a strict local minimum $x_0$ with
$\Phi(x_0)=J(x_0)=0$. Finally, setting
$$\eta=
\max \left \{ 0,
\limsup_{\|x\|\to +\infty}{{J(x)}\over {\Phi(x)}},
\limsup_{x\to x_0}{{J(x)}\over {\Phi(x)}}\right \}\ ,$$
$$\theta=\sup_{x\in \Phi^{-1}(]0,+\infty[)}{{J(x)}\over {\Phi(x)}} 
\ ,$$
assume that $\eta<\theta$.\par
Then, for each compact interval $[a,b]\subset 
\left ] {{1}\over {\theta}}, {{1}\over {\eta}}\right [$ 
(with the conventions ${{1}\over {0}}=+\infty$,
${{1}\over {+\infty}}=0$) 
there exists $r>0$ with
the following property: for every $\lambda\in [a,b]$ and every $C^1$ functional
$\Psi:X\to {\bf R}$ with compact derivative, there exists $\hat\mu>0$ such
that, for each $\mu\in [0,\hat\mu]$, the equation
$$\Phi'(x)=\lambda J'(x)+\mu\Psi'(x)$$
has at least three solutions whose norms are less than $r$.}\par
\medskip
The aim of this paper is to provide an application of Theorem A to problem $(P_{\varphi})$.\par
\smallskip
Our main result reads as follows:\par
\medskip
THEOREM 1. - {\it Let $f\in {\cal A}_p$ be such that
$$
\max\left \{ 0, \limsup_{\xi\to 0}{{\sup_{x\in \Omega}F(x,\xi)}\over
{|\xi|^p}}, \limsup_{|\xi|\to +\infty}{{\sup_{x\in \Omega}F(x,\xi)}\over
{|\xi|^p}}\right \}<$$
$$<{{\inf_{\Omega}\alpha}\over
{\int_{\Omega}\alpha(x)dx}}
\sup_{\xi\in {\bf R}\setminus\{0\}}{{\int_{\Omega}F(x,\xi)dx}\over
{|\xi|^p}}\ ,\eqno{(1)}$$
where
$$F(x,\xi)=\int_0^{\xi}f(x,s)ds$$
for all $(x,\xi)\in \Omega\times {\bf R}$. \par
Then, if we set
$$\gamma={{\int_{\Omega}\alpha(x)dx}\over {p}}\inf
\left \{
{{|\xi|^p}\over {\int_{\Omega}F(x,\xi)dx}} : \int_{\Omega}F(x,\xi)dx>0
\right \}$$
and
$$\delta={{\inf_{\Omega}\alpha}\over
p\hskip 2pt{\max\left \{ 0, \limsup_{\xi\to 0}{{\sup_{x\in \Omega}F(x,\xi)}\over
{|\xi|^p}}, \limsup_{|\xi|\to +\infty}{{\sup_{x\in \Omega}F(x,\xi)}\over
{|\xi|^p}}\right \}}}\ ,$$
for each compact interval $[a,b]\subset ]\gamma,\delta[$ there exists
$r>0$ with the following property: for every $\lambda\in [a,b]$ and
every $g\in {\cal A}_p$ there exists $\hat\mu>0$ such that, for each $\mu\in [0,\hat\mu]$, the problem
$$\cases {-\hbox {\rm div}(|\nabla u|^{p-2}\nabla u)+\alpha(x)|u|^{p-2}u=
\lambda f(x,u)+\mu g(x,u)
 & in
$\Omega$\cr & \cr {{\partial u}\over {\partial\nu}}=0 & on
$\partial \Omega$\cr} $$
has at least three weak solutions whose norms in
$W^{1,p}(\Omega)$ are less than $r$.}\par
\smallskip
PROOF.
 We are going to apply Theorem A taking $X=W^{1,p}(\Omega)$,
and
$$\Phi(u)={{1}\over {p}}\|u\|^p\ ,$$
$$J(u)=\int_{\Omega}F(x,u(x))dx$$
for all $u\in X$.\par
\smallskip
Note that, by a classical result on uniformly convex spaces,
 $\Phi\in {\cal W}_{W^{1,p}(\Omega)}$.
\smallskip
Also, set
$$\rho_1=\limsup_{\xi\to 0}{{\sup_{x\in \Omega}F(x,\xi)}\over
{|\xi|^p}}\ ,$$
$$\rho_2=
\limsup_{|\xi|\to +\infty}
{{\sup_{x\in \Omega}F(x,\xi)}\over
{|\xi|^p}}\ .$$
Fix $\epsilon>0$ and choose
 $\delta_1, \delta_2>0$ so that
$$F(x,\xi)\leq (\rho_1+\epsilon) |\xi|^p \eqno{(2)}$$
for all $(x,\xi)\in \Omega\times [-\delta_1,\delta_1]$, and
$$F(x,\xi)\leq (\rho_2+\epsilon) |\xi|^p \eqno{(3)}$$
for all $(x,\xi)\in \Omega\times ({\bf R}\setminus [-\delta_2,\delta_2])$.
Due to $(2)$ and $(3)$, since $F$ is bounded on on each bounded subset of $\Omega\times
{\bf R}$, we can choose $c>0$ 
such that
$$F(x,\xi)\leq (\rho_1+\epsilon+c)|\xi|^p$$
for all $(x,\xi)\in \Omega\times {\bf R}$. So, if we
fix $q>p$ (with $q\leq {{pn-n+p}\over {n-p}}$ if $p<n$),
by continuous
embedding, for some $c_1>0$ (independent of $\epsilon$), one has
$$J(u)\leq (\rho_1+\epsilon)\int_{\Omega}|u(x)|^pdx+c_1\|u\|^q$$
for all $u\in X$. Hence, for each $u\in X\setminus \{0\}$, we have
$${{J(u)}\over {\|u\|^p}}\leq {{(\rho_1+\epsilon)\int_{\Omega}|u(x)|^pdx}\over
{\|u\|^p}}+c_1\|u\|^{q-p}\leq {{\rho_1+\epsilon}\over {
\inf_{\Omega}\alpha}}+c_1\|u\|^{q-p}$$
and hence
$$\limsup_{u\to 0}{{J(u)}\over {\Phi(u)}}\leq
{{p(\rho_1+\epsilon)}\over {\inf_{\Omega}\alpha}}\ .\eqno{(4)}$$
Further, by $(3)$ again, for each $u\in X\setminus \{0\}$, we have
$${{J(u)}\over {\|u\|^p}}={{\int_{\Omega(|u|\leq \delta_2)}F(x,u(x))dx}\over
{\|u\|^p}}+{{\int_{\Omega(|u|>\delta_2)}F(x,u(x))dx}\over
{\|u\|^p}}\leq$$
$$\leq
{{\hbox {\rm meas}(\Omega)\sup_{\Omega\times [-\delta_2,\delta_2]}F}
\over {\|u\|^p}}+{{(\rho_2+\epsilon)\int_{\Omega}|u(x)|^pdx}\over
{\|u\|^p}}\leq
{{\hbox {\rm meas}(\Omega)\sup_{\Omega\times [-\delta_2,\delta_2]}F}
\over {\|u\|^p}}+
{{\rho_2+\epsilon}\over {\inf_{\Omega}\alpha}}$$
and hence
$$\limsup_{\|u\|\to +\infty}{{J(u)}\over {\Phi(u)}}\leq
{{p(\rho_2+\epsilon)}\over {\inf_{\Omega}\alpha}}
\ . \eqno{(5)}$$
Since $\epsilon$ is arbitrary, $(4)$ and $(5)$ give
$$\max \left \{ 
\limsup_{\|u\|\to +\infty}{{J(u)}\over {\Phi(u)}},
\limsup_{u\to 0}{{J(u)}\over {\Phi(u)}}\right \}
\leq
{{p}\over {\inf_{\Omega}\alpha}}\max\{\rho_1, \rho_2\}
\ .\eqno{(6)}$$
Moreover, it is clear that
$${{p}\over {\int_{\Omega}\alpha(x)dx}}
\sup_{\xi\in {\bf R}\setminus\{0\}}{{\int_{\Omega}F(x,\xi)dx}\over
{|\xi|^p}}\leq 
\sup_{u\in X\setminus \{0\}}{{J(u)}\over {\Phi(u)}}\ .\eqno{(7)}$$
Then, from $(1)$, $(6)$ and $(7)$, we get
$$\max \left \{0,  
\limsup_{\|u\|\to +\infty}{{J(u)}\over {\Phi(u)}},
\limsup_{u\to 0}{{J(u)}\over {\Phi(u)}}\right \}<
\sup_{u\in X\setminus \{0\}}{{J(u)}\over {\Phi(u)}}\ .$$
Hence, all the assumptions of Theorem A are satisfied (with
$x_0=0$).
The conclusion then follows taking $\Psi=J_g$ and remarking that,
due to $(6)$ and $(7)$, one has $[\gamma,\delta]\subseteq
\left [ {{1}\over {\theta}},{{1}\over {\eta}}\right ]$.
\hfill $\bigtriangleup$
\medskip
We conclude pointing out an application of Theorem 1 in a case where
$f$ is a polynomial in $\xi$.\par
\medskip
PROPOSITION 1. - {\it
Let $n=3$, let $\beta:\Omega\to {\bf R}$ be a bounded measurable function, with $\inf_{\Omega}\beta\geq 0$ and $\int_{\Omega}\beta(x)dx>0$, and let $a, b, c\in {\bf R}$, with $c>0$
and $a>-{{2b^2}\over {9c}}$. When $a>0$, assume that
$${{\sup_{\Omega}\beta\int_{\Omega}\alpha(x)dx}\over
{\inf_{\Omega}\alpha\int_{\Omega}\beta(x)dx}}<
1+{{2b^2}\over {9ac}}\ .\eqno{(8)}$$
Set
$$\gamma_0={{9c\int_{\Omega}\alpha(x)dx}\over {(9ac+2b^2)\int_{\Omega}
\beta(x)dx}}$$
and
$$\delta_0={{\inf_{\Omega}\alpha}\over
{\max \{0,a\}\sup_{\Omega}\beta}}\ .$$
Then, for each compact interval $A\subset ]\gamma_0,\delta_0[$, there
exists
$r>0$ with the following property: for every $\lambda\in A$ and
every $g\in {\cal A}_2$ there exists $\hat\mu>0$ such that, for each $\mu\in [0,\hat\mu]$, the problem
$$\cases {-\Delta u+\alpha(x)u=
\lambda \beta(x)(au+bu^2-cu^3)+\mu g(x,u)
 & in
$\Omega$\cr & \cr {{\partial u}\over {\partial\nu}}=0 & on
$\partial \Omega$\cr} $$
has at least three weak solutions whose norms in
$W^{1,2}(\Omega)$ are less than $r$.}\par
\smallskip
PROOF. Set
$$f(x,\xi)=\beta(x)(a\xi+b\xi^2-c\xi^3)$$
for all $(x,\xi)\in \Omega\times {\bf R}$. Note that $f\in {\cal A}_2$ as $n=3$. Moreover, if we set
$$\Gamma=\left \{ \xi\in {\bf R} : 
{{a}\over {2}}+{{b}\over {3}}\xi-{{c}\over {4}}\xi^2\geq 0\right \}\ $$
we clearly have
$${{\sup_{x\in \Omega}F(x,\xi)}\over {\xi^2}}=\cases
{\sup_{\Omega}\beta\left ( {{a}\over {2}}+{{b}\over {3}}\xi-{{c}\over {4}}\xi^2\right ) & if $\xi\in \Gamma$\cr & \cr
\inf_{\Omega}\beta\left (
{{a}\over {2}}+{{b}\over {3}}\xi-{{c}\over {4}}\xi^2\right ) & if
$\xi\in {\bf R}\setminus \Gamma$\ .\cr}$$
 From this, we obtain
$$\lim_{|\xi|\to +\infty}{{\sup_{x\in \Omega}F(x,\xi)}\over {\xi^2}}\leq 0$$
(since $c>0$ and $\inf_{\Omega}\beta\geq 0$) as well as
$$\lim_{\xi\to 0}{{\sup_{x\in \Omega}F(x,\xi)}\over {\xi^2}}=
\cases {{{a\sup_{\Omega}\beta}\over {2}} & if $a>0$\cr & \cr
{{a\inf_{\Omega}\beta}\over {2}} & if $a\leq 0$\ .\cr}$$
Consequently, we have
$$
\max\left \{ 0, \lim_{\xi\to 0}{{\sup_{x\in \Omega}F(x,\xi)}\over
{\xi^2}}, \lim_{|\xi|\to +\infty}{{\sup_{x\in \Omega}F(x,\xi)}\over
{\xi^2}}\right \}={{\max\{0,a\}\sup_{\Omega}\beta}\over {2}}\ .
\eqno{(9)}$$
On the other hand, we also have
$$\sup_{\xi\in {\bf R}\setminus \{0\}}{{\int_{\Omega}F(x,\xi)dx}\over
{\xi^2}}=\left ( {{a}\over {2}}+{{b^2}\over {9c}}\right ) \int_{\Omega}\beta(x)dx\ .\eqno{(10)}$$
At this point, the validity of $(1)$ (for $p=2$) follows from $(8)$, $(9)$ and $(10)$ when $a>0$. When $a\leq 0$, $(1)$ follows from
$(9)$ and $(10)$ as the right-hand side of $(10)$ is positive, by
assumption. So, $f$ satisfies the assumptions of Theorem 1, with $p=2$. Now, the conclusion follows as $\gamma_0=\gamma$ and $\delta_0=\delta$.
\hfill $\bigtriangleup$
\bigskip
\bigskip
\centerline {\bf References}\par
\bigskip
\bigskip
\noindent
[1]\hskip 5pt G. ANELLO and G. CORDARO, {\it Existence of solutions
of the Neumann problem involving the p-Laplacian via a variational
principle of Ricceri}, Arch. Math. (Basel), {\bf 79} (2002),
274-287.\par
\smallskip
\noindent
[2]\hskip 5pt G. ANELLO and G. CORDARO, {\it An existence theorem for the
Neumann problem involving the $p$-Laplacian}, J. Convex Anal., {\bf 10}
(2003), 185-198.\par
\smallskip
\noindent
[3]\hskip 5pt G. ANELLO and G. CORDARO, {\it Infinitely many positive
solutions for the Neumann problem involving the $p$-Laplacian}, Colloq.
Math., {\bf 97} (2003), 221-231.\par
\smallskip
\noindent
[4]\hskip 5pt G. ANELLO, {\it Existence of infinitely many weak solutions
for a Neumann problem}, Nonlinear Anal., {\bf 57} (2004), 199-209.\par
\smallskip
\noindent
[5]\hskip 5pt G. ANELLO, {\it A multiplicity theorem for critical points
of functionals on reflexive Banach spaces}, Arch. Math. (Basel),
{\bf 82} (2004), 172-179.\par
\smallskip
\noindent
[6]\hskip 5pt G. ANELLO, {\it Existence and multiplicity of solutions to a
perturbed Neumann problem}, Math. Nachr., {\bf 280} (2007), 1755-1764.\par
\smallskip
\noindent
[7]\hskip 5pt G. ANELLO, {\it Remarks on a multiplicity theorem for a
perturbed Neumann problem}, J. Math. Anal. Appl., {\bf 346} (2008),
274-279.\par
\smallskip
\noindent
[8]\hskip 5pt D. AVERNA and G. BONANNO, {\it Three solutions for a
Neumann boundary
value problem involving the $p$-Laplacian}, Matematiche, {\bf 60}
(2005), 81-91. \par
\smallskip
\noindent
[9]\hskip 5pt G. BILOTTA, {\it Existence of infinitely many solutions for a
quasilinear Neumann problem}, Panamer. Math. J., {\bf 13} (2003), 1-18.\par
\smallskip
\noindent
[10]\hskip 5pt G. BONANNO, {\it Multiple solutions for a Neumann boundary
value problem}, J. Nonlinear Convex Anal., {\bf 4} (2003), 287-290.\par
\smallskip
\noindent
[11]\hskip 5pt G. BONANNO, {\it Some remarks on a three critical points
theorem}, Nonlinear Anal., {\bf 54} (2003), 651-665.\par
\smallskip
\noindent
[12]\hskip 5pt G. BONANNO and P. CANDITO, {\it Three solutions to a Neumann
problem for elliptic equations involving the $p$-Laplacian}, Arch. Math.
(Basel), {\bf 80} (2003), 424-429.\par
\smallskip
\noindent
[13]\hskip 5pt F. CAMMAROTO, A. CHINN\`I and B. DI BELLA, {\it Some multiplicity
results for quasilinear Neumann problems}, Arch. Math. (Basel), {\bf 86} (2006),
154-162.\par
\smallskip
\noindent
[14]\hskip 5pt P. CANDITO, {\it Infinitely many solutions to the Neumann
problem for elliptic equations involving the $p$-Laplacian and with
discontinuous nonlinearities}, Proc. Edinb. Math. Soc., {\bf 45}
(2002), 397-409.\par
\smallskip
\noindent
[15]\hskip 5pt G. CORDARO, {\it Multiple solutions to a perturbed Neumann
problem}, Studia Math., {\bf 178} (2007), 167-175.\par
\smallskip
\noindent
[16]\hskip 5pt G. DAI, {\it Infinitely many solutions for a Neumann-type differential inclusion problem involving the $p(x)$-Laplacian},
Nonlinear Anal., to appear.\par
\smallskip
\noindent
[17]\hskip 5pt A. G. DI FALCO, {\it Infinitely many solutions of the
Neumann problem for quasilinear elliptic systems}, Matematiche, {\bf 58}
(2003), 117-130.\par
\smallskip
\noindent
[18]\hskip 5pt S. EL MANOUNI and M. KBIRI ALAOUI, {\it A result on elliptic systems with Neumann conditions via Ricceri's three critical point theorem}, Nonlinear Anal., to appear.\par
\smallskip
\noindent
[19]\hskip 5pt X. L. FAN and S. G. DENG, {\it Remarks on Ricceri's
variational principle and applications to $p(x)$-Laplacian equations},
Nonlinear Anal., {\bf 67} (2007), 3064-3075.\par
\smallskip
\noindent
[20]\hskip 5pt X. L. FAN and C. JI, {\it Existence of infinitely many
solutions for a Neumann problem involving the $p(x)$-Laplacian},
J. Math. Anal. Appl., {\bf 334} (2007), 248-260.\par
\smallskip
\noindent
[21]\hskip 5pt F. FARACI, {\it Multiplicity results for a Neumann problem
involving the $p$-Laplacian}, J. Math. Anal. Appl., {\bf 277} (2003),
180-189.\par
\smallskip
\noindent
[22]\hskip 5pt F. FARACI and A. KRIST\'ALY, {\it On an open question of Ricceri
concerning a Neumann problem}, Glasg. Math. J., {\bf 49} (2007),
189-195.\par
\smallskip
\noindent
[23]\hskip 5pt A. IANNIZZOTTO, {\it A sharp existence and localization
theorem for a Neumann problem}, Arch. Math. (Basel), {\bf 82} (2004),
352-360.\par
\smallskip
\noindent
[24]\hskip 5pt A. KRIST\'ALY, M. MIH\u{A}ILESCU and V. R\u{A}DULESCU,
{\it Two non-trivial solutions for a non-homogeneous Neumann problem:
an Orlicz-Sobolev setting}, Proc. Royal Soc. Edinburgh Sect. A, to
appear.\par
\smallskip
\noindent
[25]\hskip 5pt A. KRIST\'ALY, N. S. PAPAGEORGIOU and Cs. VARGA,
{\it Multiple solutions for a class of Neumann elliptic problems
on compact Riemannian manifolds with boundary}, Canad. Math. Bull.,
to appear.\par
\smallskip
\noindent
[26]\hskip 5pt Q. LIU, {\it Existence of three solutions for $p(x)$-Laplacian
equations}, Nonlinear Anal., {\bf 68} (2008), 2119-2127.\par
\smallskip
\noindent
[27]\hskip 5pt S. A. MARANO and D. MOTREANU, {\it Infinitely many critical
points of non-differentiable functions and applications to a Neumann type
problem involving the $p$-Laplacian}, J. Differential Equations, {\bf 182}
(2002), 108-120.\par
\smallskip
\noindent
[28]\hskip 5pt M. MIH\u{A}ILESCU, {\it Existence and multiplicity of solutions
for a Neumann problem involving the $p(x)$-Laplace operator}, Nonlinear Anal.,
{\bf 67} (2007), 1419-1425.\par
\smallskip
\noindent
[29]\hskip 5pt D. S. MOSCHETTO, {\it A quasilinear Neumann problem involving the $p(x)$-Laplacian}, Nonlinear Anal., to appear.\par
\smallskip
\noindent
[30]\hskip 5pt B. RICCERI, {\it A general variational principle and
some of its applications}, J. Comput. Appl. Math., {\bf 113}
(2000), 401-410.\par
\smallskip
\noindent
[31]\hskip 5pt B. RICCERI, {\it On a three critical points theorem},
 Arch. Math. (Basel), {\bf 75} (2000), 220-226.\par
\smallskip
\noindent
[32]\hskip 5pt B. RICCERI, {\it Infinitely many solutions of the Neumann
problem for elliptic equations involving the p-Laplacian},
Bull. London Math. Soc., {\bf 33} (2001), 331-340.\par
\smallskip
\noindent
[33]\hskip 5pt B. RICCERI, {\it Sublevel sets and global minima of coercive
functionals and local minima of their perturbations}, J. Nonlinear Convex
Anal., {\bf 5} (2004), 157-168.\par
\smallskip
\noindent
[34]\hskip 5pt B. RICCERI, {\it A multiplicity theorem for the Neumann problem},
Proc. Amer. Math. Soc., {\bf 134} (2006), 1117-1124.\par
\smallskip
\noindent
[35]\hskip 5pt B. RICCERI, {\it A three critical points theorem revisited}, Nonlinear Anal., to appear.\par
\smallskip
\noindent
[36]\hskip 5pt B. RICCERI, {\it A further three critical points theorem}, preprint.\par
\smallskip
\noindent
[37]\hskip 5pt X. SHI and X. DING, {\it Existence and multiplicity of solutions for a general $p(x)$-Laplacian Neumann problem}, Nonlinear Anal., to appear.\par
\smallskip
\noindent

\bigskip
\bigskip
\bigskip
\bigskip
Department of Mathematics\par
University of Catania\par
Viale A. Doria 6\par
95125 Catania\par
Italy\par
{\it e-mail address}: ricceri@dmi.unict.it

\bye

\bye

\bye

\bye